# SOBRE UN MÉTODO DE REGULARIZACIÓN PARA IDENTIFICAR UNA FUENTE EN UNA ECUACIÓN ELÍPTICA

Guillermo Federico Umbricht[&], Diana Rubio [*]

Centro de Matemática Aplicada, ECyT

Universidad Nacional de San Martin

[&]guilleungs@yahoo.com.ar   [*]diarubio@gmail.com

**Resumen**

El objetivo de este trabajo es estudiar numéricamente el desempeño de un método de regularización. Esta técnica fue desarrollada para resolver el problema mal planteado de la estimación de una fuente unidimensional para una ecuación de Poisson de dos dimensiones a partir de mediciones tomadas a lo largo de una curva del dominio. El método consiste en introducir en la ecuación un término de regularización que depende de un parámetro, denominado parámetro de regularización. Aquí se presentan resultados para distintos valores de este parámetro como así también para diferentes niveles de ruido en los datos utilizados para la estimación y se discuten algunas consideraciones sobre su aplicación efectiva.

**Palabras Clave**: Ecuación de Poisson, Problema mal planteado, Método de regularización, Estimación de la fuente.

# IDENTIFICATION OF THE SOURCE IN A POISSON EQUATION USING A METHOD OF REGULARIZATION

**Abstract**

The aim of this paper is to numerically study the performance of a method of regularization. This technique was developed to solve the ill-posed problem of estimating a source-dimensional Poisson equation for two dimensions from measurements taken over a line inside the domain. The proposed method consists in adding a regularization term to the equation that depends on a parameter, which is called regularization parameter. In this paper we show the results for different values of this parameter as well as for different levels of noise in the data used for estimation. After analyzing the results some considerations on its effective implementation are discussed.

### *Elección del Tema*

El problema inverso de la identificación de la fuente tiene aplicaciones directas en diversas disciplinas como economía, física, ingenierías, química, biología, salud, entre otros. Los desarrollos y resultados matemáticos pueden luego ser adaptados a problemas particulares de aplicación de las distintas ciencias. Por lo tanto, el interés no es sólo académico, sino de aplicación concreta para la resolución de problemas reales.

### *Definición del Problema*

Estudiamos numéricamente la performance del método de regularización introducido por Xiao-Xiao Li et al. en [12] para resolver el problema inverso que consiste en estimar la fuente (desconocida) $f(x)$ de la ecuación de Poisson con condición inicial y de frontera:

$$\begin{aligned} -u_{xx} - u_{yy} &= f(x), \quad -\infty < x < +\infty, \ 0 < y < +\infty \\ u(x,0) &= 0, -\infty < x < +\infty \\ u(x,y)|_{y \to \infty} &\ acotada, \quad -\infty < x < +\infty \\ u(x,1) &= g(x), \quad -\infty < x < +\infty, \end{aligned} \quad (1)$$

donde $f(x)$ sólo depende de una variable espacial y $g(x)$ representa los valores a ser medidos. En general, los datos de medidos no son exactos, la precisión depende en la mayoría de los casos, del aparato utilizado para las mediciones. Estos errores de medición se los denomina *ruido* y se deben tener en cuenta en la resolución del problema. Consideramos $g_\delta$ como la función de las mediciones (ruidosas) obtenidas, que satisface:

$$||g - g_\delta||_{L^2(\mathbb{R})} \leq \delta \quad (2)$$

Aquí la constante $\delta$ representa el nivel de ruido de los datos de entrada.

### ***Antecedentes (Marco Teórico o Conceptual)***

La teoría del problema inverso fue inicialmente desarrollada por investigadores que centran su estudio en la Geofísica y buscaban entender los fenómenos que ocurren en el interior de la tierra basados en un conjunto de datos obtenidos en la superficie.

Los métodos de regularización surgieron a partir de la necesidad de encontrar soluciones apropiadas, estables, y la imposibilidad de obtenerlas a partir problema que se quiere resolver. Entre los métodos clásicos se encuentran la regularización de Tikhonov, la regresión Ridge, la descomposición en valores singulares truncada (DVST), entre muchos otros.

### ***Justificación del Estudio***

El problema de la identificación de la fuente desconocida ha sido muy estudiado y ha recibido considerable atención de muchas investigaciones actuales debido esencialmente a que se han encontrado disímiles y múltiples aplicaciones en los distintos campos de la ciencia tales como la conducción de calor [10], la identificación de fisuras [13], la teoría electromagnética [5], la prospección geofísica [3], la detección de contaminantes [7], la tomografía y la tomografía de impedancia eléctrica [1] y determinación de tumores en un tejido biológico [8], entre otros.

La ecuación de Poisson es una ecuación diferencial en derivadas parciales con un amplio uso en electrostática [2], ingeniería mecánica [9] y en Física teórica [4]. El problema inverso de estimar la fuente en la ecuación de Poisson es un problema *mal planteado* o *mal condicionado* en el sentido de Hadamard [6], ya que las soluciones obtenidas no son estables, esto es, la solución no varía de forma continua con los datos del modelo. La problemática que genera la falta de estabilidad es que la solución hallada puede encontrarse muy distante de la verdadera. Por este motivo el problema se resuelve utilizando un *método de*

*regularización* que consiste en el agregado de un término en la ecuación a estudiar [12].

Una característica común de las técnicas de regularización es que todas ellas requieren la utilización de un parámetro que depende del problema a regularizar, por lo general se denomina *parámetro de regularización.* Este parámetro controla el peso dado a la minimización del término penalizador.

Nos interesa analizar la performance del método de regularización mencionado bajo distintas situaciones

### *Alcances del trabajo*

Para este trabajo nos concentramos en un problema de identificación de una fuente imponiendo características de suavidad de la misma para una ecuación de tipo Poisson con condiciones iniciales y de frontera dadas.

### *Objetivos*

- Analizar numéricamente el comportamiento de un método de regularización.
- Discutir variantes del método estudiado con el fin de mejorar su efectividad.

**Estimación de la fuente en una ecuación de tipo Poisson**

Consideramos la ecuación de Poisson:

$$-u_{xx} - u_{yy} = f(x), \quad -\infty < x < +\infty, \quad 0 < y < +\infty$$

con las siguientes condiciones inicial y de frontera

$$u(x, 0) = 0, -\infty < x < +\infty$$

$$u(x,y)|_{y\to\infty}\ acotada, \quad -\infty < x < +\infty$$
$$u(x,1) = g(x), \quad -\infty < x < +\infty,$$

donde $f(x)$ es la fuente desconocida que sólo depende de una variable espacial y $g(x)$ (supuestamente conocida o con valores conocidos) representa los valores reales correspondientes a las mediciones.

El problema consiste en hallar la fuente $f(x)$ conociendo únicamente los datos obtenidos a partir de los valores de la función $u$ en los puntos donde $y = 1$, en otras palabras, hallar $f(x)$ conociendo $g(x)$. La solución del problema está dado por (ver [11], [12])

$$f(x) = \frac{1}{\sqrt{2\pi}} \int_{-\infty}^{\infty} e^{i\xi x} \frac{\xi^2}{(1-e^{-\xi})} \hat{g}(\xi) d\xi \qquad (3)$$

donde $\hat{g}$ denota la transformada de Fourier de la función $g$.

La presencia del factor $\xi^2/(1-e^{-\xi})$ hace inestable la solución, ya que este factor aumenta cuando se incrementa $\xi$, resultando en un factor de amplificación de $g(x)$, y por lo tanto motivo la solución $f(x)$ no es continua con respecto a los datos, por lo que el problema está mal planteado. Esto hace imposible encontrar una solución confiable y precisa al problema con métodos clásicos por lo que requiere el uso de alguna técnica de regularización.

### El método de regularización

El método de regularización propuesto en [12] consiste en resolver la ecuación que se obtiene al agregar un término con la derivada de segundo orden de la fuente $f(x)$ que se quiere estimar, con las mismas condiciones iniciales y de frontera, es decir

$$-u_{xx} - u_{yy} + \mu^2 f_{xx}(x) = f(x), \quad -\infty < x < +\infty, \ 0 < y < +\infty$$
$$u(x,0) = 0, -\infty < x < +\infty \qquad (4)$$

$$u(x,y)|_{y\to\infty} \text{ acotada}, \quad -\infty < x < +\infty$$
$$u(x,1) = g_\delta(x), \quad -\infty < x < +\infty,$$

donde el parámetro $\mu$ es el parámetro de regularización y corresponde al peso o grado de importancia del término $f_{xx}$ y $g_\delta(x)$ representa los datos u observaciones, con errores de medición de nivel $\delta$.

La ecuación puede ser reformulada en el espacio de frecuencias utilizando la Transformada de Fourier (ver [12]). Luego antitransformando se obtiene la solución regularizada:

$$f_{\delta,\mu}(x) = \frac{1}{\sqrt{2\pi}} \int_{-\infty}^{\infty} e^{i\xi x} \frac{\xi^2}{(1-e^{-\xi})(1+\xi^2\mu^2)} \hat{g}_\delta(\xi) d\xi \tag{5}$$

Los subíndices $\delta, \mu$ en la función $f$ dada en (5) indican la dependencia de la solución con respecto al nivel de ruido de los datos $\delta$ y del parámetro de regularización $\mu$.

En [12] se propone elegir el parámetro de regularización en base a la relación entre el nivel de ruido de los datos y una cota de la fuente a estimar para obtener una cota del error de la aproximación, según se exhibe en el siguiente teorema.

**Teorema** [12] Sea $f(x)$ la solución exacta del problema (1), $f_{\delta,\mu}(x)$ la aproximación regularizada definida en (5) y $g_\delta(x)$ los datos medidos en $y=1$. Se supone que la $f(x)$ y $g_\delta(x)$ satisfacen (4). Se define

$$\mu = \left(\frac{\delta}{E}\right)^{1/(p+2)} \tag{6}$$

donde $p > 0$ y $E > 0$ es una constante real que verifica

$$\|f(\cdot)\|_{H^p(\mathbb{R})} := \left(\int_{-\infty}^{\infty} |\hat{f}(\xi)|^2 (1+\xi^2)^p d\xi\right)^{1/2} \leq E. \tag{7}$$

Entonces se obtiene la siguiente estimación para el error de aproximación

$$\|f(\cdot) - f_{\delta,\mu}(\cdot)\| \leq 2\delta^{p/(p+2)} \left(1 + \frac{1}{2} max\{1, \mu^{2-p}\}\right) \tag{8}$$

### *__Resultados__*

Comentamos aquí algunos aspectos generales de la resolución numérica. El lector interesado en la implementación del método puede ver los detalles en [11].

Con el fin de estudiar la bondad de la aproximación de la fuente, construimos un ejemplo con funciones $u(x,y)$, $f(x)$ que satisfacen el problema (1). De este modo, al conocer la función a estimar, se podrán calcular los errores de estimación y analizar la bondad del método.

Considerando
$$u(x,y) = (1 - e^{-y})cos(x)$$
se puede ver que satisface la ecuación y condiciones dadas en (1) para la $f(x) = cos(x)$, es decir

$$-u_{xx} - u_{yy} = cos(x), \quad -\infty < x < +\infty, \quad 0 < y < +\infty$$
$$u(x,0) = 0, \quad -\infty < x < +\infty,$$
$$u(x,y)|_{y\to\infty} \; acotada, \quad -\infty < x < +\infty.$$

Los datos son $N$ mediciones ruidosas de la función $u$ en $(x,1)$, de donde
$$g_\delta(x_j) = (1 - e^{-1})\cos(x_j) + \varepsilon(\delta) \qquad j = 1,..,N$$
siendo $\varepsilon(\delta)$ una variable aleatoria con distribución normal, media 0 y desvío estándar $\delta$ (nivel de ruido). Con la finalidad de obtener $\hat{g}_\delta(\xi)$ se utiliza la función *fft* para calcular la transformada discreta de Fourier de $g_\delta$.

Se utilizó el software MATLAB para la resolución numérica de los ejemplos. Los datos $g_\delta(x)$, se simularon a partir de los valores de la función $g(x) = u(x,1)$ sumando una perturbación distribuida normalmente, con media 0 y desvio estándar $\delta$. Esta perturbación simula los errores cometidos al tomar mediciones experimentales y fue generada por la función *randn*.

La solución aproximada $f_{\delta,\mu}(x)$ es calculada utilizando la fórmula dada en (5), es decir

$$f_{\delta,\mu}(x) = \frac{1}{\sqrt{2\pi}} \int_{-\infty}^{\infty} e^{i\xi x} \frac{\xi^2}{(1-e^{-\xi})(1+\xi^2\mu^2)} \hat{g}_\delta(\xi) d\xi$$

Luego calculamos el error de aproximación, $||f_{\delta,\mu}(x) - f(x)||$, considerando diferentes niveles de ruido.

Con respecto a los valores del parámetro $\mu$, se usaron los propuestos en [12] usando la fórmula (6) con $E = 1$, de donde $\mu = \delta^{1/(p+2)}$.

Se consideraron distintos valores para el parámetro de regularización $\mu$ y se graficaron las soluciones aproximadas obtenidas.
En la Figura 1 se muestra la función $f(x) = cos(x)$, que se quiere estimar, junto con la solución obtenida sin regularizar ($\mu = 0$) calculada mediante la fórmula (3):

$$f(x) = \frac{1}{\sqrt{2\pi}} \int_{-\infty}^{\infty} e^{i\xi x} \frac{\xi^2}{(1-e^{-\xi})} \hat{g}_\delta(\xi) d\xi.$$

Se puede observar que la solución encontrada sin regularizar tiene la forma de la solución pero presenta una gran cantidad de oscilaciones mientras que las soluciones regularizadas obtenidas con la fórmula dada en (5):

$$f_{\delta,\mu}(x) = \frac{1}{\sqrt{2\pi}} \int_{-\infty}^{\infty} e^{i\xi x} \frac{\xi^2}{(1-e^{-\xi})(1+\xi^2\mu^2)} \hat{g}_\delta(\xi) d\xi,$$

resultan más suaves, como se puede observar en las Figuras 2-4.

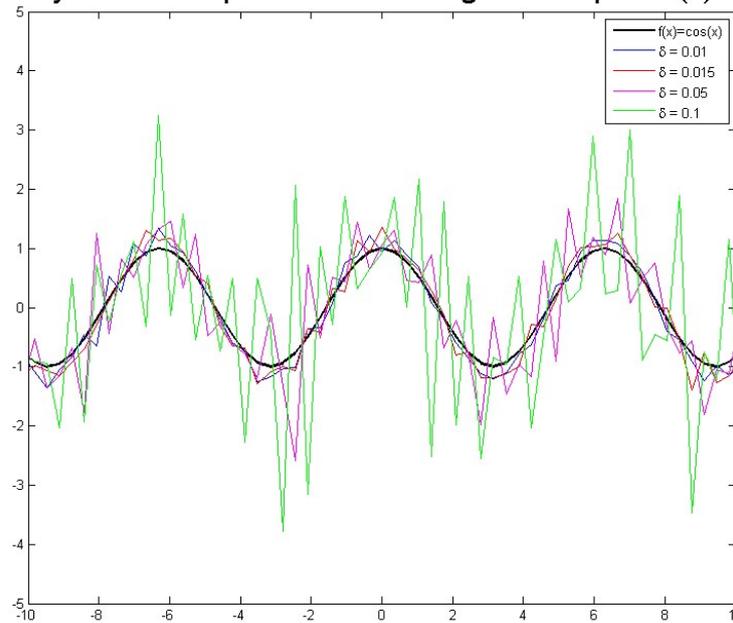

**Figura 1.** Cálculo de la fuente $f(x) = cos(x)$ sin usar regularización, obtenidas a partir de datos tomados de $g_\delta(x) = (1 - e^{-1})cos(x) + \varepsilon(\delta)$ para distintos niveles de ruido $\delta$.

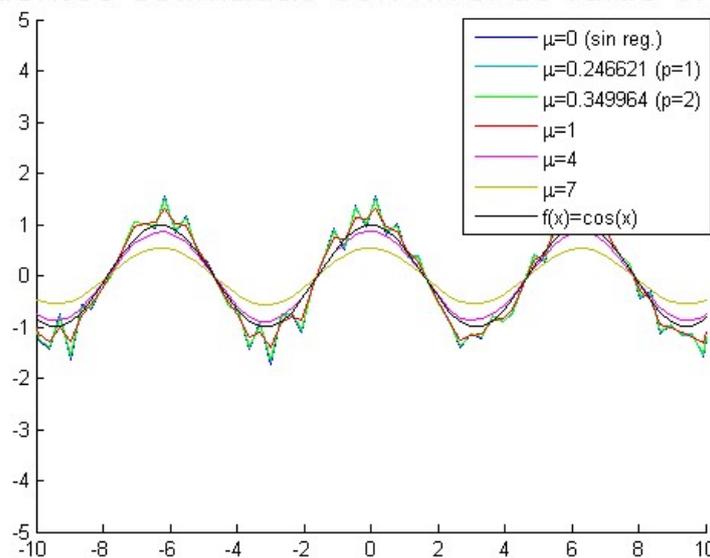

**Figura 2.** Soluciones aproximadas para $f(x) = cos(x)$ obtenidas a partir de datos tomados de $g_\delta(x) = (1 - e^{-1})cos(x) + \varepsilon(\delta)$ para distintos parámetros de regularización con nivel de ruido $\delta = 0.015$.

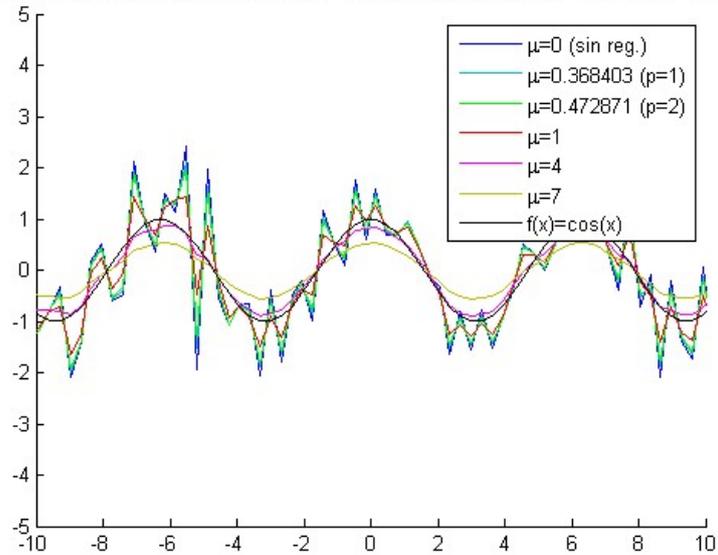

**Figura 3.** Soluciones aproximadas para $f(x) = cos(x)$ obtenidas a partir de datos tomados de $g_\delta(x) = (1 - e^{-1})cos(x) + \varepsilon(\delta)$ para distintos parámetros de regularización con nivel de ruido $\delta = 0.05$.

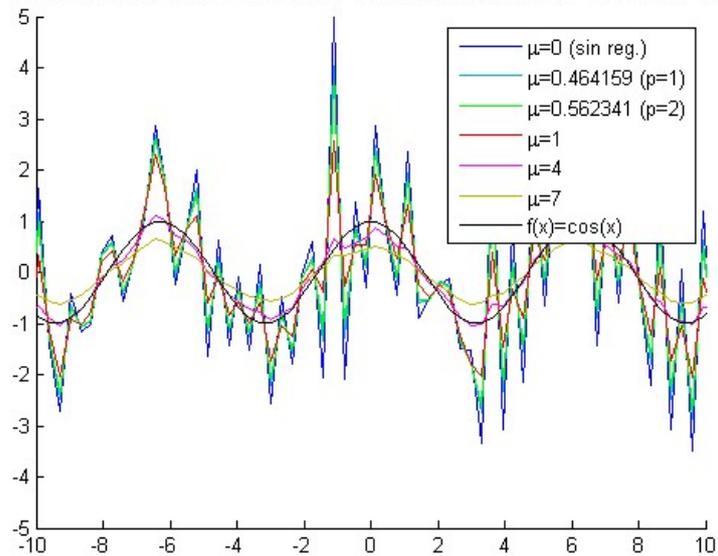

**Figura 4.** Soluciones aproximadas para $f(x) = cos(x)$ obtenidas a partir de datos tomados de $g_\delta(x) = (1 - e^{-1})cos(x) + \varepsilon(\delta)$ para distintos parámetros de regularización con nivel de ruido $\delta = 0.1$.

## *Discusiones*

En el método estudiado los autores proponen considerar como parámetro de regularización el dado en la ecuación (6) con $E = 1$, esto es

$$\mu = \delta^{1/(p+2)} \qquad (9)$$

con $p > 0$. Los niveles $\delta$ de ruido son acotados, es decir que se encuentran dentro de un cierto rango de valores. Si se hace variar $p$ para cada posible valor $\delta$ se tiene que $\mu^2 = \delta$ para $p = 0$ y $\mu^2 \sim 1$ para $p \to \infty$ (será más parecido a 1 cuanto mayor sea el valor de $p$). Luego, para cada valor de $\delta$ resulta

$$\delta \leq \mu^2 \leq 1. \qquad (10)$$

En las simulaciones presentadas, además de considerar esta regla de selección del parámetro para $p = 1$ y $p = 2$, se consideraron valores de $\mu$ mayores a 1 obteniéndose aproximaciones más precisas. Observamos que estos valores no pueden ser obtenidos a partir de la fórmula dada para el parámetro, pero nada impide que puedan ser utilizados para la regularización. Recordemos que el parámetro de regularización interviene en la fórmula para el cálculo de la solución regularizada

$$f_{\delta,\mu}(x) = \frac{1}{\sqrt{2\pi}} \int_{-\infty}^{\infty} e^{i\xi x} \frac{\xi^2}{(1-e^{-\xi})(1+\xi^2\mu^2)} \hat{g}_\delta(\xi) d\xi$$

y teóricamente representa el peso del término $f_{xx}(x)$ en la ecuación (4) que se debe resolver.

Ampliamos ahora este análisis considerando valores de $\mu$ entre 0 y 40 ($0 \leq \mu 2 \leq 1600$) muy superiores a los valores propuestos en [12]. En la Figura 5 mostramos los resultados de los errores obtenidos de la solución regularizada al estimar $f(x) = \cos(x)$ con distintos niveles de ruido de observación.

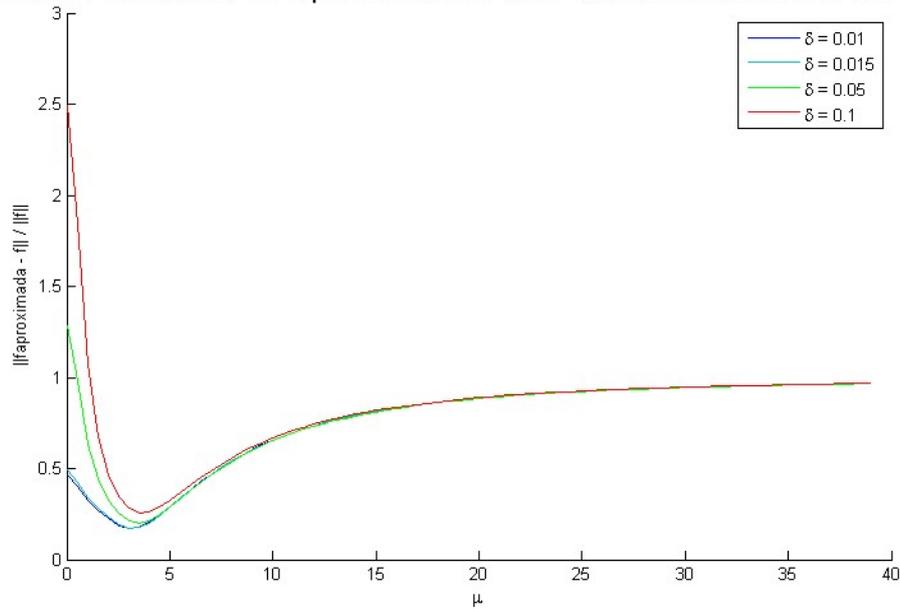

**Figura 5.** Errores de aproximación de $f(x) = \cos(x)$ estimadas con el método de regularización en función del parámetro $\mu$.

Observando los resultados obtenidos, se puede apreciar que el método de regularización estudiado provee una estimación más suavizada de la obtenida sin regularizar. Por otro lado, al utilizar los valores de parámetro de regularización propuesto por sus autores (ver [12]), se puede observar una mejora. Sin embargo la mejora puede ser aún mayor si se consideran valores mayores para el parámetro de regularización, según se pudo observar en los gráficos presentados aquí. Más aún, en la última figura se puede apreciar que los errores de aproximación hallados al comienzo disminuyen rápidamente con el valor del parámetro $\mu$, mostrando un valor óptimo que no varía demasiado con el nivel de ruido. Las curvas de error sugieren que en este ejemplo se puede obtener una buena aproximación eligiendo un valor cercano a 3 para el parámetro $\mu$.

Por otro lado, ya que el método consiste en modificar la ecuación introduciendo un término de la forma $\mu^2 f_{xx}(x)$, el mismo supone que la fuente es suficientemente suave. Además, si la derivada segunda se anula, el problema no se modifica. Por este motivo pareciera carecer de sentido tratar de usar el método para estimar funciones lineales a trozos, es decir funciones que no son derivables en todos los puntos y además, se anulan en las regiones donde es derivable, sin embargo ejemplos de este tipo fueron incluidos en [12].

## ***Conclusiones***

En este trabajo consideramos la identificación de una fuente desconocida (que es función de una variable) en la ecuación de Poisson de dos dimensiones y observamos que se trata de un problema mal planteado, es decir, que en caso de existir solución, ésta no depende de manera continua de los datos de entrada. Consideramos el método de regularización propuesto en [12] a fin de obtener una solución aproximada.

Analizamos numéricamente las soluciones obtenidas para distintos valores del parámetro de regularización teniendo en cuenta distintos niveles de ruido de medición de los datos.

Se requiere un estudio más profundo para poder dar conclusiones acerca del mejor valor del parámetro, que seguramente, va a depender de la función a estimar y de los datos utilizados para la estimación.

## ***Bibliografía***